\documentclass{amsart}
\usepackage{amssymb,amsmath,amsthm}
\usepackage[utf8]{inputenc}
\usepackage[italian,english]{babel}
\usepackage{mathrsfs}
\usepackage[textmath,displaymath,floats,graphics,tightpage]{preview}
\usepackage[all]{xy}
\usepackage{pgf,tikz}
\usetikzlibrary{arrows}
\usepackage{float}
\usepackage[labelformat=empty]{caption,subfig}

\restylefloat{figure}

\theoremstyle{definition}
\newtheorem{defin}{Definition}[section]
\theoremstyle{plain}
\newtheorem{thm}[defin]{Theorem}

\newtheorem{prop}[defin]{Proposition}

\newtheorem{step}{Step}
\theoremstyle{remark}
\newtheorem{ex}[defin]{Example}
\newtheorem{remark}[defin]{Remark}

\newcounter{num}
\newcommand{\dom}{\addtocounter{num}{1}\thenum)}

\title{On the Hilbert function on $\mathbb P^1\times \mathbb P^1$}

\author{Paola Bonacini}
\email{bonacini@dmi.unict.it}

\author{Lucia Marino}
\email{lmarino@dmi.unict.it}
\address{Università degli Studi di Catania\\
  Viale A. Doria 6\\set
95125 Catania\\
Italy}

\begin{document}

\maketitle

\begin{abstract}
Let $Q=\mathbb P^1\times \mathbb P^1$ and let
$X\subset Q$ be a $0$-dimensional scheme. This paper is a
first step towards the characterization of Hilbert functions of
$0$-dimensional schemes in $Q$. In particular we show how,
under some conditions on $X$, its Hilbert function changes when we add points to $X$ lying on
a $(1,0)$ or $(0,1)$-line. As a particular case we show also that if $X$ is ACM
this result holds without any additional hypothesis.
\end{abstract}

\section{Introduction}

Let $Q=\mathbb P^1\times \mathbb P^1$ and let
$X\subset Q$ be a $0$-dimensional scheme. Let $R$ and $C$ be,
respectively, a $(1,0)$ and a $(0,1)$-line not containing any point of $X$ and let
$Z$ be a $0$-dimensional scheme given by $X$ and some points on $R$ or
$C$. In this paper we deal with the problem of finding the Hilbert matrix
(function) of $Z$ with respect to the Hilbert matrix of $X$. A first approach
was given in a very particular case in 1992 in \cite{GMR}, with the only 
perspective of comprehending the Hilbert functions of ACM
$0$-dimensional schemes in $Q$. This paper is the first real step towards the characterization of
the Hilbert functions of $0$-dimensional schemes in $Q$ that  are not ACM.

In Theorem \ref{T:1} and Theorem \ref{T:2} we improve
the result in \cite{GMR}, under some geometric and algebraic conditions that, as we see in
Example \ref{ex}, can not be suppressed without any further
assumption. In Theorem
\ref{T:5} and Theorem \ref{T:6} we see
that the result holds for all ACM schemes $X$ without any additional
condition. As an application in Section \ref{Esempio} we compute the
Hilbert matrix of any non ACM reduced set of points in $Q$ having a certain position in a grid of $(1,0)$ and $(0,1)$-lines. This
previously could be done and was known just for ACM $0$-dimensional schemes.  
  
A good reference for a general discussion on $0$-dimensional schemes on $\mathbb P^1\times
\mathbb P^1$ is \cite{GMR}, in which there are the most important
results about the Hilbert function. Further results on the Hilbert
function has been obtained just in the particular case of fat
points (see for example \cite{GVT} and \cite{VT2}).  

\section{Notation and preliminary results}

Let $k$ be an algebraically closed field,let $\mathbb P^1=\mathbb
P^1_k$, let $Q=\mathbb P^1\times \mathbb P^1$ and let
$\mathscr O_Q$ be its structure sheaf. For any sheaf $\mathscr F$ we
denote:
\[
\mathscr F(a,b)=\mathscr F\otimes \mathscr O_Q(a,b).
\]
Let us consider the bi-graded ring:
\[
S=H^0_*(\mathscr O_Q)=\bigoplus_{a,b\ge 0}H^0(\mathscr O_Q(a,b)).
\]
For any bi-graded $S$-module $N$ let $N_{(i,j)}$ the component of
degree $(i,j)$.

Given $X\subset Q$ $0$-dimensional scheme, let $I(X)\subset S$ be the
associated saturated sheaf, $S(X)=S/I(X)$ and $\mathscr I_X\subset
\mathscr O_Q$ its ideal sheaf. 
\begin{defin}
  The function:
\[
M_X\colon \mathbb Z\times \mathbb Z\rightarrow \mathbb N
\]
defined by:
\[
M_X(i,j)=\dim_k {S(X)}_{(i,j)}=\dim_k S_{(i,j)}-\dim_k {I(X)}_{(i,j)}
\]
is called \emph{Hilbert function} of $X$. The function $M_X$ can be
represented as a matrix with infinite integers entries: 
\[
M_X=(M_X(i,j))=(m_{ij})
\]
called \emph{Hilbert matrix} of $X$. 
\end{defin}
Note that $M_X(i,j)=0$ for either $i<0$ or $j<0$ and so we restrict
ourselves to the range $i\ge 0$ and $j\ge 0$. Moreover, for $i\gg 0$
and $j\gg 0$ $M_X(i,j)=\deg X$.

\begin{defin}
 Given the Hilbert matrix $M_X$ of a $0$-dimensional scheme $X\subset Q$, the \emph{first difference of the Hilbert function} of $X$ is the
  matrix:
\[
\Delta M_X=(c_{ij}),
\]
where:
\[
c_{ij}=m_{ij}-m_{i-1j}-m_{ij-1}+m_{i-1j-1}.
\]
\end{defin}

We consider the following matrices:
\[
\Delta^R M_X=(a_{ij})\mbox{ and }\Delta^C M_X=(b_{ij}),
\]
with $a_{ij}=m_{ij}-m_{ij-1}$ and $b_{ij}=m_{ij}-m_{i-1j}$. Note that:
\[
c_{ij}=a_{ij}-a_{i-1j}=b_{ij}-b_{ij-1}
\]
and
\[
m_{ij}=\sum_{\substack{h\le i\\ k\le j}}c_{hk}.
\]

\begin{thm}[{\cite[Theorem 2.11]{GMR}}]  \label{T0}
  Given a $0$-dimensional scheme $X\subset Q$ and given its Hilbert matrix
  $M_X$, the first difference $\Delta M_X=(c_{ij})$ satisfies the
  following conditions:
  \begin{enumerate}
  \item $c_{ij}\le 1$ and $c_{ij}=0$ for $i\gg 0$ or $j\gg 0$;
\item if $c_{ij}\le 0$, then $c_{rs}\le 0$ for any $(r,s)\ge (i,j)$;
\item for every $(i,j)$ $0\le \sum_{t=0}^j c_{it}\le
  \sum_{t=0}^jc_{i-1t}$ and $0\le \sum_{t=0}^i c_{tj}\le \sum_{t=0}^i c_{tj-1}$.
  \end{enumerate}
\end{thm}

\begin{remark}  \label{rm}
If $X\subset Q$ is a $0$-dimensional scheme, let us consider $a=\min\{i\in \mathbb
N\mid I(X)_{(i,0)}\ne 0\}-1$ and $b=\min\{j\in \mathbb
N\mid I(X)_{(0,j)}\ne 0\}-1$. Then by Theorem \ref{T0} $\Delta M_X$ is zero out of the
rectangle with opposite vertices $(0,0)$ and $(a,b)$. In this case we
say that $\Delta M_X$ is of syze $(a,b)$.
\end{remark}

  Let $M_X=(m_{ij})$ be the Hilbert matrix of a $0$-dimensional scheme
  $X\subset Q$. Using the notation in \cite{GMR}, for every $j\ge 0$ we set:
\[
i(j)=\min\{t\in \mathbb N\mid m_{tj}=m_{t+1j}\}=\min\{t\in \mathbb
N\mid b_{t+1j}=0\}.
\]
and for every $i\ge 0$ we set:
\[
j(i)=\min\{t\in \mathbb N\mid m_{it}=m_{it+1}\}=\min\{t\in \mathbb
N\mid a_{it+1}=0\}
\]
In particular, we see that $i(0)=a$ and $j(0)=b$. 

Let $X\subset Q$ be a $0$-dimensional scheme and let $L$ be a line
defined by a form $l$. Let $J=(I(X),l)$ and let $d=\deg(\operatorname{sat} J)$. Then we call $d$ the number of points
of $X$ on the line $L$ and, by abuse of notation, we make the position
$d=\#(X\cap L)$. We say that $L$ is disjoint from $X$ if $d=0$.

The key result used in Section \ref{main} is the following:
\begin{thm}[{\cite[Theorem 2.12]{GMR}}] \label{T}
  Let $X\subset Q$ be a $0$-dimensional scheme and let $M_X=(m_{ij})$
  be its Hilbert matrix. Then for every $j\ge 0$ there are just
  $a_{i(0)j}-a_{i(0)j+1}$ lines of type $(1,0)$ each containing just
  $j+1$ points of $X$ and, similarly, for every $i\ge 0$ there are
  just $b_{ij(0)}-b_{i+1j(0)}$ lines of type $(0,1)$ each containing
  just $i+1$ points of $X$.
\end{thm}

The result that in this paper we improve is given by the following:

\begin{thm}[{\cite[Lemma 2.15]{GMR}}]  \label{T2}
 Let $X\subset Q$ be a 0-dimensional scheme and let $M_X$ be its Hilbert
matrix. Let $R_0$,\dots, $R_a$
and $C_0$,\dots,$C_b$ be, respectively, the
$(1,0)$ and $(0,1)$-lines containing $X$ and at least one point of $X$. 
 Let $R$ be a $(1,0)$-line disjoint from $X$ and let $Y=X\cup
R\cap(C_0\cup\dots\cup C_n)$, with $n\ge b$ and
$C_{b+1}$,\dots,$C_n$ arbitrary $(0,1)$-lines. Then:
\[
\Delta M_Y^{(i,j)}=
\begin{cases}
 1 & \text{for } i=0,\, j\le n\\
0 & \text{for }i=0,\, j\ge n+1\\
\Delta M_X^{(i-1,j)} & \text{for } i\ge 1. 
\end{cases}
\]
\end{thm}
Of course a similar result can be proved by adding $m+1$ points on a
$(0,1)$-line $C$ disjoint from $X$. So, with the previous notation, it
is possible to prove the following result. 

\begin{thm}
 Let $C$ be a $(0,1)$-line disjoint from $X$. Let $Y=X\cup
C\cap(R_0\cup\dots\cup R_m)$, $m\ge a$, and
$R_{a+1}$,\dots,$R_m$ arbitrary $(1,0)$-lines. Then:
\[
\Delta M_Y^{(i,j)}=
\begin{cases}
 1 & \text{for } i\le m,\,j=0\\
0 & \text{for } i\ge m+1,\, j=0\\
\Delta M_X^{(i,j-1)} & \text{for } j\ge 1. 
\end{cases}
\]
\end{thm}

\section{The first difference of the Hilbert function}  \label{main}

Let $X\subset Q$ be a 0-dimensional scheme and let $M_X$ be its Hilbert
matrix. In all this paper we suppose that $\Delta M_X$ is of size
$(a,b)$ and we denote by $R_0$,\dots, $R_a$
and $C_0$,\dots,$C_b$, respectively, the
$(1,0)$ and $(0,1)$-lines containing $X$ and at least one point of $X$. 

\begin{thm}  \label{T:1}
Let $R$ be a $(1,0)$-line disjoint from $X$. Let $C_{b+1}$,\dots,$C_n$,
$n\ge b$, be arbitrary $(0,1)$-lines and $i_1$,\dots,$i_r\in
\{0,\dots,b\}$. Let $\mathcal P=\{R\cap C_i\mid i\in\{0,\dots,n\},\,
i\ne i_1,\dots,i_r\}$ and let $Z=X\cup \mathcal P$. Suppose also that on the $(0,1)$-line $C_{i_k}$
there are $q_k$ points of $X$ for $k=1,\dots,r$ and that $q_1\le
q_2\le \dots \le q_r$. Then, given
$T=\{(q_1,n),(q_2,n-1),\dots,(q_r,n-r+1)\}$, we have:
\[
\Delta M_Z^{(i,j)}=
\begin{cases}
1 & \text{if } i=0,\, j\le n\\
0 & \text{if }i=0,\, j\ge n+1\\
 \Delta M_X^{(i-1,j)} & \text{if } i\ge 1\text{ and }(i,j)\notin T\\
 \Delta M_X^{(i-1,j)}-1 & \text{if } i\ge 1\text{ and }(i,j)\in T
\end{cases}
\]
if one of the following conditions holds:
\begin{enumerate}
\item $r=1$;
\item $r\ge 2$, $q_{r-1}<q_r$ and for any $k\in \{1,\dots,r-1\}$ and $i\ge q_k$\, $\Delta
  M_X^{(i,n-k+1)}=0$;
\item $r\ge 2$, $q_{r-1}=q_r$ and for any $k\in \{1,\dots,r\}$ and $i\ge q_k$\, $\Delta
  M_X^{(i,n-k+1)}=0$.
\end{enumerate}
\end{thm}
\begin{proof}
Let $Y=X\cup (R\cap (\bigcup_{i=0}^n C_i))$. By Theorem \ref{T2} it is
sufficient to prove that:
\[
\Delta M_Z(i,j)=
\begin{cases}
  \Delta M_Y(i,j) & \text{if }(i,j)\notin T\\
\Delta M_Y(i,j)-1& \text{if }(i,j)\in T.
\end{cases}
\]
We divide the proof in different steps.

\begin{step}  \label{s:1}
   $\Delta M_Z^{(0,j)}=\Delta M_Y^{(0,j)}=1$ for $j\le n$, $\Delta
   M_Z^{(0,j)}=\Delta M_Y^{(0,j)}=0$ for $j\ge n+1$ and $\Delta M_Z^{(i,j)}=\Delta
M_Y^{(i,j)}=\Delta M_X^{(i-1,j)}$ for any $(i,j)$ with $j<n-r+1$ and $i\ge 1$.
\end{step}
It is easy to see that $\Delta M_Z^{(0,j)}=1$ for $j\le n$, because
for such values of $j$ $h^0(\mathscr I_Z(0,j))=0$. Moreover, $\Delta
M_Z^{(0,j)}=0$ for $j\ge n+1$ by Remark \ref{rm}.

Taken $(i,j)$, with $j< n-r+1$ and $i\ge 1$, any $(i,j)$-curve containing $Z$ must
contain $R$ and so $h^0(\mathscr
I_Z(i,j))=h^0(\mathscr I_X(i-1,j))$ and $\Delta M_Z^{(i,j)}=\Delta
M_X^{(i-1,j)}$.
\\

Let $r_1$,\dots,$r_{t+1}$ be a sequence of positive integers such
that
$q_{1}=\dots=q_{r_1}<q_{r_1+1}=\dots=q_{r_2}<\dots<q_{r_t+1}=\dots=q_{r_{t+1}}=q_r$
and let $r_0=0$.  
 
\begin{step}  \label{s:2}
If $h\in \{1,\dots,t+1\}$, then $\Delta M_Z^{(i,j)}=\Delta
M_Y^{(i,j)}$ for any $(i,j)\le (q_{r_{h-1}+1}-1,n-r_{h-1})$ and for
$(i,j)=(q_{r_{h-1}+1},j)$ with $j<n-r_h+1$.
\end{step}
Taken $(i,j)\le (q_{r_{h-1}+1}-1,n-r_{h-1})$, then any $(i,j)$-curve
containing $Z$ must contain $C_{i_{r_{h-1}+1}}$,\dots,$C_{i_r}$ and so it must contain
$R$. This means that $h^0(\mathscr I_Z(i,j))=h^0(\mathscr I_Y(i,j))$
and so that $\Delta M_Z^{(i,j)}=\Delta M_Y^{(i-1,j)}$.

Taken $(i,j)=(q_{r_{h-1}+1},j)$ with $j<n-r_h+1$, then any
$(i,j)$-curve containing $Z$ must contain
$C_{i_{r_h+1}}$,\dots,$C_{i_r}$ and so it must contain $R$. Again this
implies $h^0(\mathscr
I_Z(i,j))=h^0(\mathscr I_Y(i,j))$ and $\Delta M_Z^{(i,j)}=\Delta
M_Y^{(i,j)}$.

\begin{step}  \label{s:3}
For any $1\le h\le t+1$ one of the following conditions holds:
\begin{enumerate}
\item there exists $\overline j$ with $n-r_h+1\le \overline j\le n-r_{h-1}$ such that  $\Delta
M_Z^{(q_{r_{h-1}+1},j)}=\Delta M_Y^{(q_{r_{h-1}+1},j)}$ for any $j< \overline j$ and
$\Delta M_Z^{(q_{r_{h-1}+1},\overline j)}<\Delta M_Y^{(q_{r_{h-1}+1},\overline j)}$;
\item  $\Delta
M_Z^{(q_{r_{h-1}+1},j)}=\Delta M_Y^{(q_{r_{h-1}+1},j)}$ for any $n-r_h+1\le j\le n-r_{h-1}$.
\end{enumerate}
\end{step}

Since $Z\subset Y$ we see that $M_Z(q_{r_{h-1}+1},n-r_h+1)\le
M_Y(q_{r_{h-1}+1},n-r_h+1)$. Moreover, by Step \ref{s:2} we
see  that $M_Z(i,j)=M_Y(i,j)$ for any $(i,j)<(q_{r_{h-1}+1},n-r_h+1)$. This
implies that:
\[
\Delta M_Z^{(q_{r_{h-1}+1},n-r_h+1)}\le \Delta M_Y^{(q_{r_{h-1}+1},n-r_h+1)}.
\]
If $\Delta M_Z^{(q_{r_{h-1}+1},n-r_h+1)}=\Delta M_Y^{(q_{r_{h-1}+1},n-r_h+1)}$, then we
can repeat the previous procedure to show that $\Delta
M_Z^{(q_{r_{h-1}+1},n-r_h+2)}\le \Delta M_Y^{(q_{r_{h-1}+1},n-r_h+2)}$. By iterating this
procedure we get the conclusion of Step \ref{s:3}.

\begin{step} \label{s:4}
\
\begin{enumerate}
\item If $r\ge 2$ and $q_{r-1}=q_r$, given $h\in \{1,\dots,t+1\}$ and $j\in\{n-r_h+1,\dots,n-r_{h-1}\}$,
  we have:
\[
\sum_{i=q_{r_{h-1}+1}}^{a+1}\Delta M_Z^{(i,j)}=\Delta M_Y^{(q_{r_{h-1}+1},j)}-1;
\]
\item If $r\ge 2$ and $q_{r-1}<q_r$, given $h\in \{1,\dots,t\}$ and $j\in\{n-r_h+1,\dots,n-r_{h-1}\}$,
  we have:
\[
\sum_{i=q_{r_{h-1}+1}}^{a+1}\Delta M_Z^{(i,j)}=\Delta M_Y^{(q_{r_{h-1}+1},j)}-1
\]
and 
\[
\sum_{i=q_r}^{a+1}\Delta M_Z^{(i,n-r+1)}=\sum_{i=q_r}^{a+1}\Delta M_Y^{(i,n-r+1)}-1.
\]
\end{enumerate}
\end{step}

Let us first note that by Theorem \ref{T}: 
\[
  a_{i(0)n-r}(Z)-a_{i(0)n-r+1}(Z)=
\sum_{i\le a+1}
\Delta M_Z^{(i,n-r)}-\sum_{i\le a+1} \Delta M_Z^{(i,n-r+1)}
\]
is equal to the number of $(1,0)$-lines containing precisely $n-r+1$ points of
$Z$, while: 
\[  
a_{i(0)n-r}(Y)-a_{i(0)n-r+1}(Y)=
\sum_{i\le a+1}
\Delta M_Y^{(i,n-r)}-\sum_{i\le a+1} \Delta M_Y^{(i,n-r+1)}
\]
is equal to the number of $(1,0)$-lines containing precisely $n-r+1$ points of
$Y$. By hypothesis it must be:
\begin{multline*}
  \sum_{i\le a+1}
\Delta M_Z^{(i,n-r)}-\sum_{i\le a+1} \Delta M_Z^{(i,n-r+1)}=\\
=\sum_{i\le a+1}
\Delta M_Y^{(i,n-r)}-\sum_{i\le a+1} \Delta M_Y^{(i,n-r+1)}+1.
\end{multline*}
By Step \ref{s:1} this implies that:
\[
\sum_{i\le a+1} \Delta M_Z^{(i,n-r+1)}=\sum_{i\le a+1} \Delta M_Y^{(i,n-r+1)}-1.
\]
Let us now suppose that for some $j\ge n-r+1$, with $j<n$, we have:
\begin{equation}
  \label{eq:4}
  \sum_{i\le a+1} \Delta M_Z^{(i,j)}=\sum_{i\le a+1} \Delta M_Y^{(i,j)}-1.
\end{equation}
We will show that:
\begin{equation}
  \label{eq:5}
  \sum_{i\le a+1} \Delta M_Z^{(i,j+1)}=\sum_{i\le a+1} \Delta M_Y^{(i,j+1)}-1.
\end{equation}
Again, by Theorem \ref{T} $\sum_{i\le a+1}
\Delta M_Z^{(i,j)}-\sum_{i\le a+1} \Delta M_Z^{(i,j+1)}$ is equal
to the number of $(1,0)$-lines containing precisely $j+1$ points of
$Z$, while $\sum_{i\le a+1}
\Delta M_Y^{(i,j)}-\sum_{i\le a+1} \Delta M_Y^{(i,j+1)}$ is equal
to the number of $(1,0)$-lines containing precisely $j+1$ points of
$Y$. By hypothesis it must be:
\[
  \sum_{i\le a+1}
\Delta M_Z^{(i,j)}-\sum_{i\le a+1} \Delta M_Z^{(i,j+1)}=
\sum_{i\le a+1}
\Delta M_Y^{(i,j)}-\sum_{i\le a+1} \Delta M_Y^{(i,j+1)}.
\]
By \eqref{eq:4} it means that \eqref{eq:5} holds, so that:
\[
\sum_{i\le a+1} \Delta M_Z^{(i,j)}=\sum_{i\le a+1} \Delta M_Y^{(i,j)}-1
\]
for any $j$ with $n-r+1\le j\le n$. Now the hypotheses on $X$, Step \ref{s:1} and Step \ref{s:2} give us the conclusion of Step \ref{s:4}.

\begin{step} \label{s:5}
 If $r\ge 2$, $\Delta M_Z^{(i,j)}=\Delta
M_Y^{(i,j)}=0$ for any $(i,j)\ge (q_1+1,n-r_1+1)$ and $\Delta M_Z^{(i,j)}=\Delta
M_Y^{(i,j)}-1$ 
for any $(i,j)\in \{(q_1,n),(q_1,n-1),\dots,(q_1,n-r_1+1)\}$.
\end{step}

By Theorem \ref{T} we know that:
\[
b_{q_1-1j(0)}(Z)-b_{q_1j(0)}(Z)=\sum_{j\le n} \Delta
M_Z^{(q_1-1,j)}-\sum_{j\le n} \Delta M_Z^{(q_1,j)}
\] 
is equal to the number
of $(0,1)$-lines containing exactly $q_1$ points of $Z$ and, in the
same way that:
\[
b_{q_1-1j(0)}(Y)-b_{q_1j(0)}(Y)=\sum_{j\le n} \Delta
M_Y^{(q_1-1,j)}-\sum_{j\le n} \Delta M_Y^{(q_1,j)}
\] 
is equal to the number of $(0,1)$-lines containing exactly $q_1$ points of
$Y$. So by construction we have:
\begin{equation}   \label{eq:1}
\sum_{j\le n} \Delta
M_Z^{(q_1-1,j)}-\sum_{j\le n} \Delta M_Z^{(q_1,j)}=\sum_{j\le n} \Delta
M_Y^{(q_1-1,j)}-\sum_{j\le n} \Delta M_Y^{(q_1,j)}+r_1.
\end{equation}
By what we proved in Step \ref{s:2} we see that:
\[
\sum_{j\le n} \Delta
M_Z^{(q_1-1,j)}=\sum_{j\le n} \Delta
M_Y^{(q_1-1,j)}
\]
so that by \eqref{eq:1} and again by Step \ref{s:2}
we get:
\begin{equation}
  \label{eq:2}
  \sum_{n-r_1+1}^{n} \Delta M_Z^{(q_1,j)}=\sum_{n-r_1+1}^{n} \Delta M_Y^{(q_1,j)}-r_1.
\end{equation}

By this equality and by Step \ref{s:3} we see that there exists $n-r_1+1\le \overline j\le
n$ such that $\Delta
M_Z^{(q_1,j)}=\Delta M_Y^{(q_1,j)}$ for any $j<\overline j$ and
$\Delta M_Z^{(q_1,\overline j)}<\Delta M_Y^{(q_1,\overline j)}$. In
particular, $\Delta M_Z^{(q_1,\overline j)}\le 0$ and so by Theorem \ref{T0} we have:
\[
\Delta M_Z^{(i,\overline j)}\le 0
\]
for any $i\ge q_1$. By Step \ref{s:4} we have:
\[
\sum_{i=q_1}^{a+1} \Delta M_Z^{(i,\overline j)}=\Delta
M_Y^{(q_1,\overline j)}-1
\]
\[
\Rightarrow 0\ge \sum_{i=q_1+1}^{a+1} \Delta M_Z^{(i,\overline
  j)}=\Delta M_Y^{(q_1,\overline j)}-\Delta M_Z^{(q_1,\overline
  j)}-1\ge 0.
\]
This means that $\Delta M_Y^{(q_1,\overline j)}-\Delta M_Z^{(q_1,\overline
  j)}-1=0$ and that $\Delta M_Z^{(i,\overline j)}=0$ for any $i\ge
q_1+1$. 

Now take any $j>\overline j$. By the fact that $\Delta M_Z^{(i,\overline j)}=0$ for any $i\ge
q_1+1$ and by Theorem \ref{T0} we can say that $\Delta
M_Z^{(i,j)}\le 0$ for any $i\ge q_1+1$ and any $j>\overline j$. By
Step \ref{s:4} we get:
\begin{equation}
  \label{eq:10}
  0\ge \sum_{i=q_1+1}^{a+1} \Delta M_Z^{(i,j)}=\Delta M_Y^{(q_1,j)}-\Delta M_Z^{(q_1,j)}-1
\end{equation}
so that:
\[
\Delta M_Z^{(q_1,j)}\ge \Delta M_Y^{(q_1,j)}-1
\]
for any $j>\overline j$. So we can say that:
\begin{equation}
  \label{eq:6}
  \sum_{j=n-r_1+1}^n \Delta M_Z^{(q_1,j)}\ge \sum_{j=n-r_1+1}^n \Delta
M_Y^{(q_1,j)}-n+\overline j-1.
\end{equation}
This fact compared to \eqref{eq:2} gives us that $\overline j\le
n-r_1+1$, but by hypothesis $\overline j\ge n-r_1+1$ and so it must be
$\overline j=n-r_1+1$. This implies that the inequality in
\eqref{eq:6} is an equality, which means that:
\[
\Delta M_Z^{(q_1,j)}=\Delta M_Y^{(q_1,j)}-1
\]
for any $j\ge \overline j=n-r_1+1$, with $j\le n$, and by \eqref{eq:10}  $\Delta M_Z^{(i,j)}=\Delta
M_Y^{(i,j)}=0$ for any $(i,j)\ge (q_1+1,n-r_1+1)$. 

\begin{step}  \label{s:6} 
If $r\ge 2$, $q_{r-1}=q_r$ and for any $k\in \{1,\dots,r\}$ and $i\ge q_k$\,
$\Delta M_X^{(i,n-k+1)}=0$, then  $\Delta M_Z^{(i,j)}=\Delta
M_Y^{(i,j)}=0$ for any $(i,j)\ge (q_k+1,n-k+1)$ and $k\in \{1,\dots,r\}$ 
and $\Delta M_Z^{(i,j)}=\Delta
M_Y^{(i,j)}-1$ for any $(i,j)\in \{(q_1,n),(q_2,n-1),\dots,(q_{r},n-r+1)\}$.
\end{step}

We proceed iterating the procedure given in Step \ref{s:5}. So let us suppose that for some $h\in \{2,\dots,t+1\}$
the equalities in the claim hold for any $i<q_{r_{h-1}+1}$ and for
any $j\ge n-r_{h-1}+1$. We will show that they hold also for
$q_{r_{h-1}+1}\le i<q_{r_h+1}$ and for any $j\ge n-r_h+1$.

To this end, we repeat what we did in Step \ref{s:5}. So, as done
before, we see that $\sum_{j\le n} \Delta
M_Z^{(q_{r_{h-1}+1}-1,j)}-\sum_{j\le n} \Delta M_Z^{(q_{r_{h-1}+1},j)}$ is the number of
$(0,1)$-lines containing precisely $q_{r_{h-1}+1}$ points of $Z$, while $\sum_{j\le n} \Delta
M_Y^{(q_{r_{h-1}+1}-1,j)}-\sum_{j\le n} \Delta M_Y^{(q_{r_{h-1}+1},j)}$ is the number of
$(0,1)$-lines containing precisely $q_{r_{h-1}+1}$ points of $Y$. By hypothesis it must
be:
\begin{itemize}
\item[\dom] if $q_{r_{h-1}+1}-1>q_{r_{h-1}}$:
\begin{multline*}
\sum_{j\le n} \Delta
M_Z^{(q_{r_{h-1}+1}-1,j)}-\sum_{j\le n} \Delta
M_Z^{(q_{r_{h-1}+1},j)}=\\
=\sum_{j\le n} \Delta
M_Y^{(q_{r_{h-1}+1}-1,j)}-\sum_{j\le n} \Delta M_Y^{(q_{r_{h-1}+1},j)}+r_h-r_{h-1};
\end{multline*}
\item[\dom] if $q_{r_{h-1}+1}-1=q_{r_{h-1}}$:
\begin{multline*}
\sum_{j\le n} \Delta
M_Z^{(q_{r_{h-1}+1}-1,j)}-\sum_{j\le n} \Delta
M_Z^{(q_{r_{h-1}+1},j)}=\\
=\sum_{j\le n} \Delta
M_Y^{(q_{r_{h-1}+1}-1,j)}-\sum_{j\le n} \Delta M_Y^{(q_{r_{h-1}+1},j)}+r_h-r_{h-1}-(r_{h-1}-r_{h-2}).
\end{multline*}
\end{itemize}

By what we proved in Step \ref{s:1} and Step \ref{s:2} and by inductive
hypothesis we see that these equalities are both equivalent to the following:
\begin{equation}
  \label{eq:7}
   \sum_{n-r_h+1}^{n-r_{h-1}} \Delta
   M_Z^{(q_{r_{h-1}+1},j)}=\sum_{n-r_h+1}^{n-r_{h-1}} \Delta M_Y^{(q_{r_{h-1}+1},j)}-r_h+r_{h-1}.
\end{equation}

By this equality and by Step \ref{s:3} we see that there exists $n-r_h+1\le \overline j\le
n-r_{h-1}$ such that $\Delta
M_Z^{(q_{r_{h-1}+1},j)}=\Delta M_Y^{(q_{r_{h-1}+1},j)}$ for any $j<\overline j$ and
$\Delta M_Z^{(q_{r_{h-1}+1},\overline j)}<\Delta
M_Y^{(q_{r_{h-1}+1},\overline j)}$. In
particular, $\Delta M_Z^{(q_{r_{h-1}+1},\overline j)}\le 0$ and so by
Theorem \ref{T0} we have:
\[
\Delta M_Z^{(i,\overline j)}\le 0
\]
for any $i\ge q_{r_{h-1}+1}$. By Step \ref{s:4} we have:
\[
\sum_{i=q_{r_{h-1}+1}}^{a+1} \Delta M_Z^{(i,\overline j)}=\Delta
M_Y^{(q_{r_{h-1}+1},\overline j)}-1
\]
\[
\Rightarrow 0\ge \sum_{i=q_{r_{h-1}+1}+1}^{a+1} \Delta M_Z^{(i,\overline
  j)}=\Delta M_Y^{(q_{r_{h-1}+1},\overline j)}-\Delta M_Z^{(q_{r_{h-1}+1},\overline
  j)}-1\ge 0.
\]
This means that $\Delta M_Y^{(q_{r_{h-1}+1},\overline j)}-\Delta M_Z^{(q_{r_{h-1}+1},\overline
  j)}-1=0$ and that $\Delta M_Z^{(i,\overline j)}=0$ for any $i\ge
q_{r_{h-1}+1}+1$. 

Now take any $j>\overline j$. By the fact that $\Delta M_Z^{(i,\overline j)}=0$ for any $i\ge
q_{r_{h-1}+1}+1$ and by Theorem \ref{T0} we can say that $\Delta
M_Z^{(i,j)}\le 0$ for any $i\ge q_{r_{h-1}+1}+1$ and any $j>\overline j$. By
Step \ref{s:4} we get:
\begin{equation}
  \label{eq:11}
  0\ge \sum_{i=q_{r_{h-1}+1}+1}^{a+1} \Delta M_Z^{(i,j)}=\Delta M_Y^{(q_{r_{h-1}+1},j)}-\Delta M_Z^{(q_{r_{h-1}+1},j)}-1
\end{equation}
so that:
\[
\Delta M_Z^{(q_{r_{h-1}+1},j)}\ge \Delta M_Y^{(q_{r_{h-1}+1},j)}-1
\]
for any $j>\overline j$. So we can say that:
\begin{equation}
  \label{eq:9}
  \sum_{j=n-r_h+1}^{n-r_{h-1}} \Delta M_Z^{(q_{r_{h-1}+1},j)}\ge \sum_{j=n-r_h+1}^{n-r_{h-1}} \Delta
M_Y^{(q_{r_{h-1}+1},j)}-n+r_{h-1}+\overline j-1.
\end{equation}
This fact compared to \eqref{eq:7} gives us that $\overline j\le
n-r_h+1$, but by hypothesis $\overline j\ge n-r_h+1$ and so it must be
$\overline j=n-r_h+1$. This implies that the inequality in
\eqref{eq:9} is an equality, which means that:
\[
\Delta M_Z^{(q_{r_{h-1}+1},j)}=\Delta M_Y^{(q_{r_{h-1}+1},j)}-1
\]
for any $j\ge \overline j=n-r_h+1$, with $j\le n-r_{h-1}$, and by
\eqref{eq:11} and by inductive hypothesis
$\Delta M_Z^{(i,j)}=\Delta M_Y^{(i,j)}=0$ for any $i\ge
q_{r_{h-1}+1}$ and any $j\ge n-r_h+1$.   

In this way we have proved the conclusion holds for any $(i,j)$, with
$i<q_{r_h+1}$ and the proof works by iteration.

\begin{step} \label{s:7}
  If either $r=1$ or $r\ge 2$, $q_{r-1}<q_r$ and for any $k\in \{1,\dots,r-1\}$ and $i\ge q_k$ $\Delta
  M_X^{(i,n-k+1)}=0$, then  $\Delta M_Z^{(i,j)}=\Delta
M_Y^{(i,j)}=0$ for any $(i,j)>(q_k,n-k+1)$ and $k\in \{1,\dots,r\}$ and $\Delta M_Z^{(i,j)}=\Delta
M_Y^{(i,j)}-1$ for any $(i,j)\in \{(q_1,n),(q_2,n-1),\dots,(q_{r},n-r+1)\}$.
\end{step}

Let us first suppose that $r\ge 2$ and that $q_{r-1}<q_r$. In this
case the procedure given in Step \ref{s:6} can be repeated for any $h\in
\{2,\dots,t\}$. This means that the equalities in the conclusion of Step
\ref{s:7} hold for
any $i<q_r$, for any $j\ge n-r+2$ and also for any $j\le n-r$ by Step
\ref{s:1}, i.e. for any $j\ne n-r+1$. 

If $r=1$, then by Step \ref{s:1} and Remark \ref{rm} we see that the
conclusion holds for any $(i,j)$ with $j\ne n-r+1$ and for $(i,n)$ with
$i<q_1$. So in both cases we will show that $\Delta M_Z^{(q_r,n-r+1)}=\Delta
M_Y^{(q_r,n-r+1)}-1$ and that $\Delta M_Z^{(i,n-r+1)}=\Delta
M_Y^{(i,n-r+1)}=0$ for $i>q_r$.

As done before, we see that  $\sum_{j\le n} \Delta
M_Z^{(q_r-1,j)}-\sum_{j\le n} \Delta M_Z^{(q_r,j)}$ is the number of
$(0,1)$-lines containing precisely $q_r$ points of $Z$, while $\sum_{j\le n} \Delta
M_Y^{(q_r-1,j)}-\sum_{j\le n} \Delta M_Y^{(q_r,j)}$ is the number of
$(0,1)$-lines containing precisely $q_r$ points of $Y$. By hypothesis it must
be:
\[
\sum_{j\le n} \Delta
M_Z^{(q_r-1,j)}-\sum_{j\le n} \Delta M_Z^{(q_r,j)}=\sum_{j\le n} \Delta
M_Y^{(q_r-1,j)}-\sum_{j\le n} \Delta M_Y^{(q_r,j)}+1.
\]
Since the equalities in the claim hold for any $i<q_r$ and for any $j\ne n-r+1$, we see that this equality is equivalent to the following:
\begin{equation}
  \label{eq:8}
   \Delta M_Z^{(q_r,n-r+1)}=\Delta M_Y^{(q_r,n-r+1)}-1.
\end{equation}

Since the $(0,1)$-lines containing exactly $q_r+1$ points of $Z$ are
one less than those containing exactly $q_r+1$ points of $Y$, we see
that:
\[
\sum_{j\le n} \Delta
M_Z^{(q_r,j)}-\sum_{j\le n} \Delta M_Z^{(q_r+1,j)}=\sum_{j\le n} \Delta
M_Y^{(q_r,j)}-\sum_{j\le n} \Delta M_Y^{(q_r+1,j)}-1,
\]
which, by our hypotheses, implies:
\[
\Delta M_Z^{(q_r+1,n-r+1)}=\Delta M_Y^{(q_r+1,n-r+1)}.
\]

By iterating the procedure, taken any $i\ge q_r+2$, the $(0,1)$-lines containing exactly $i$ points of $Z$ are
also those containing exactly $i$ points of $Y$, so that:
\[
\sum_{j\le n} \Delta
M_Z^{(i-1,j)}-\sum_{j\le n} \Delta M_Z^{(i,j)}=\sum_{j\le n} \Delta
M_Y^{(i,j)}-\sum_{j\le n} \Delta M_Y^{(i,j)},
\]
which, by our hypotheses, implies:
\[
\Delta M_Z^{(i,n-r+1)}=\Delta M_Y^{(i,n-r+1)}.
\]
\end{proof}

In the same way, with the above notation, we can prove the following theorem:

\begin{thm}  \label{T:2}
Let $C$ be a $(0,1)$-line disjoint from $X$. Let $R_{a+1}$,\dots,$R_m$,
$m\ge a$, be arbitrary $(1,0)$-lines and $j_1$,\dots,$j_r\in
\{0,\dots,a\}$. Let $\mathcal P=\{C\cap R_j\mid j\in\{0,\dots,m\},\,
j\ne j_1,\dots,j_r\}$ and let $Z=X\cup \mathcal P$. Suppose also that on the $(1,0)$-line $R_{j_k}$
there are $p_k$ points of $X$ for $k=1,\dots,r$ and that $p_1\le
p_2\le \dots \le p_r$. Then, given
$T=\{(m,p_1),(m-1,p_2),\dots,(m-r+1,p_r)\}$, we have:
\[
\Delta M_Z^{(i,j)}=
\begin{cases}
1 & \text{if } i\le m,\, j=0\\
0 & \text{if }i\ge m+1,\, j=0\\
 \Delta M_X^{(i,j-1)} & \text{if } j\ge 1\text{ and }(i,j)\notin T\\
 \Delta M_X^{(i,j-1)}-1 & \text{if } j\ge 1\text{ and }(i,j)\in T
\end{cases}
\]
if one of the following conditions holds:
\begin{enumerate}
\item $r=1$;
\item $r\ge 2$, $p_{r-1}<p_r$ and for any $k\in \{1,\dots,r-1\}$ and $j\ge p_k$\, $\Delta M_X^{(m-k+1,j)}=0$;
\item $r\ge 2$, $p_{r-1}=p_r$ and for any $k\in \{1,\dots,r\}$ and $j\ge p_k$\, $\Delta M_X^{(m-k+1,j)}=0$.
\end{enumerate}
\end{thm}
\begin{proof}
  The proof works as in Theorem \ref{T:1}. 
\end{proof}

Under the notation of Theorem \ref{T:1} we prove the following:
\begin{thm}  \label{T:3}
 If one the following conditions holds:
 \begin{enumerate}
 \item $q_{r-1}<q_r$ and $n\ge b+r-1$,
\item $q_{r-1}=q_r$ and $n\ge b+r$,
 \end{enumerate}
then:
 \[
\Delta M_Z^{(i,j)}=
\begin{cases}
1 & \text{if } i=0,\, j\le n\\
0 & \text{if }i=0,\, j\ge n+1\\
 \Delta M_X^{(i-1,j)} & \text{if } i\ge 1\text{ and }(i,j)\notin T\\
 \Delta M_X^{(i-1,j)}-1 & \text{if } i\ge 1\text{ and }(i,j)\in T.
\end{cases}
\]
\end{thm}
\begin{proof}
\
\begin{enumerate}
\item By Remark \ref{rm} our hypothesis imply that $\Delta M_X^{(i,j)}=0$ for any
  $i\ge 0$ and for any $j\ge n-r+2$, so that the hypothesis of Theorem
  \ref{T:1} holds;
\item in this case by hypothesis we have that $\Delta M_X^{(i,j)}=0$ for any
  $i\ge 0$ and for any $j\ge n-r+1$, so that the hypothesis of Theorem
  \ref{T:1} holds.
\end{enumerate}
\end{proof}

In the same way under the notation of Theorem \ref{T:2} we prove the
following result:
\begin{thm}  \label{T:4}
If one the following
 conditions holds:
 \begin{enumerate}
 \item $p_{r-1}<p_r$ and $m\ge a+r-1$,
\item $p_{r-1}=p_r$ and $m\ge a+r$,
 \end{enumerate}
then:
\[
\Delta M_Z^{(i,j)}=
\begin{cases}
1 & \text{if } i\le m,\, j=0\\
0 & \text{if }i\ge m+1,\, j=0\\
 \Delta M_X^{(i,j-1)} & \text{if } j\ge 1\text{ and }(i,j)\notin T\\
 \Delta M_X^{(i,j-1)}-1 & \text{if } j\ge 1\text{ and }(i,j)\in T.
\end{cases}
\]
\end{thm}
\begin{proof}
The works as in Theorem \ref{T:3}.  
\end{proof}

\begin{ex}  \label{ex}
In these examples we will show that if the hypothesis of Theorem \ref{T:1}
does not hold, then the conclusion is not necessarily true. As a notation,
we represent the $(1,0)$-lines as horizontal lines and the
$(0,1)$-lines as vertical lines.
\begin{enumerate}
\item Let us consider a scheme $X$ union of three generic points
and its first difference $\Delta M_X$.
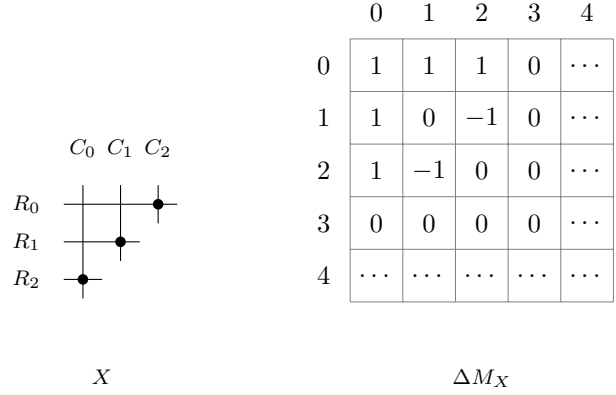
\begin{figure}[H]
\begin{preview}
\begin{center}
\subfloat[$X$]{
\begin{tikzpicture}[line cap=round,line join=round,>=triangle 45,x=0.5cm,y=0.5cm]
\clip(-1,-1) rectangle (4,6);
\fill [color=black] (1,1) circle (2pt);
\fill [color=black] (2,2) circle (2pt);
\fill [color=black] (3,3) circle (2pt);
\draw (1,0.5) -- (1,3.5);
\draw (2,1.5) -- (2,3.5);
\draw (3,2.5) -- (3,3.5);
\draw (0.5,1) -- (1.5,1);
\draw (0.5,2) -- (2.5,2);
\draw (0.5,3) -- (3.5,3);
\draw[color=black] (1,4.5) node {\footnotesize $C_0$};
\draw[color=black] (2,4.5) node {\footnotesize $C_1$};
\draw[color=black] (3,4.5) node {\footnotesize $C_2$};
\draw[color=black] (-0.5,1) node {\footnotesize $R_2$};
\draw[color=black] (-0.5,2) node {\footnotesize $R_1$};
\draw[color=black] (-0.5,3) node {\footnotesize $R_0$};
\end{tikzpicture}
}
\hspace{1cm}
\subfloat[$\Delta M_X$]{
\begin{tikzpicture}[x=0.7cm,y=0.7cm]
\clip(0,0) rectangle (7,7);
  \draw[style=help lines,xstep=1,ystep=1] (1,1) grid (6,6);
\foreach \x in {1,...,5} \draw (\x,1) +(.5,.5)  node {\dots};
\foreach \y in {2,...,5} \draw (5,\y) +(.5,.5) node {\dots};
\foreach \x in {1,...,4} \draw (\x,2.5) +(.5,0) node {$0$};
\foreach \y in {3,...,5} \draw (4.5,\y) +(0,.5) node {$0$};
\draw (1.5,3.5) node {$1$};
\draw (2.5,3.5) node {$-1$};
\draw (3.5,3.5) node {$0$};
\draw (1.5,4.5) node {$1$};
\draw (2.5,4.5) node {$0$};
\draw (3.5,4.5) node {$-1$};
\foreach \x in {1,...,3} \draw (\x,5.5) +(.5,0) node {$1$};
\foreach \x in {0,...,4} \draw (\x,6.5) +(1.5,0) node {$\x$};
\foreach \y in {0,...,4} \draw (0.5,5.5-\y) node {$\y$};
\end{tikzpicture}
}
\end{center}
\end{preview}
\end{figure} 
Let $R$ be a $(1,0)$-line disjoint from $X$ and let $Z$ be the following scheme. 
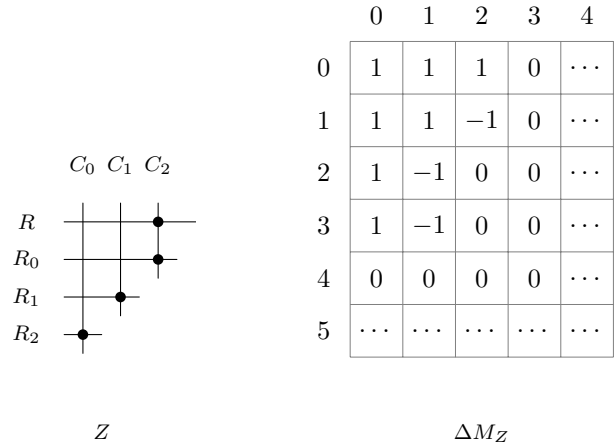
\begin{figure}[H]
\begin{preview}
\begin{center}
\subfloat[$Z$]{
\begin{tikzpicture}[line cap=round,line join=round,>=triangle 45,x=0.5cm,y=0.5cm]
\clip(-1,-1) rectangle (4,6);
\fill [color=black] (1,1) circle (2pt);
\fill [color=black] (2,2) circle (2pt);
\fill [color=black] (3,3) circle (2pt);
\fill [color=black] (3,4) circle (2pt);
\draw (1,0.5) -- (1,4.5);
\draw (2,1.5) -- (2,4.5);
\draw (3,2.5) -- (3,4.5);
\draw (0.5,1) -- (1.5,1);
\draw (0.5,2) -- (2.5,2);
\draw (0.5,3) -- (3.5,3);
\draw (0.5,4) -- (4.5,4);
\draw[color=black] (1,5.5) node {\footnotesize $C_0$};
\draw[color=black] (2,5.5) node {\footnotesize $C_1$};
\draw[color=black] (3,5.5) node {\footnotesize $C_2$};
\draw[color=black] (-0.5,1) node {\footnotesize $R_2$};
\draw[color=black] (-0.5,2) node {\footnotesize $R_1$};
\draw[color=black] (-0.5,3) node {\footnotesize $R_0$};
\draw[color=black] (-0.5,4) node {\footnotesize $R$};
\end{tikzpicture}
}
\hspace{1cm}
\subfloat[$\Delta M_Z$]{
\begin{tikzpicture}[x=0.7cm,y=0.7cm]
\clip(0,0) rectangle (7,8);
\draw[style=help lines,xstep=1,ystep=1] (1,1) grid (6,7);
\foreach \x in {1,...,5} \draw (\x,1) +(.5,.5)  node {\dots};
\foreach \y in {2,...,6} \draw (5,\y) +(.5,.5) node {\dots};
\foreach \x in {1,...,4} \draw (\x,2.5) +(.5,0) node {$0$};
\foreach \y in {3,...,6} \draw (4.5,\y) +(0,.5) node {$0$};
\draw (1.5,3.5) node {$1$};
\draw (2.5,3.5) node {$-1$};
\draw (3.5,3.5) node {$0$};
\draw (1.5,4.5) node {$1$};
\draw (2.5,4.5) node {$-1$};
\draw (3.5,4.5) node {$0$};
\foreach \x in {1,...,2} \draw (\x,5.5) +(.5,0) node {$1$};
\draw (3.5,5.5) node {$-1$};
\foreach \x in {1,...,3} \draw (\x,6.5) +(.5,0) node {$1$};
\foreach \x in {0,...,4} \draw (\x,7.5) +(1.5,0) node {$\x$};
\foreach \y in {0,...,5} \draw (0.5,6.5-\y) node {$\y$};
\end{tikzpicture}
}
\end{center}
\end{preview}
\end{figure} 
In this case, under the notation of
Theorem \ref{T:1} we have $r=2$, $n=2$, $q_1=q_2=1$ and $\Delta
M_X(q_1,n)=\Delta M_X(1,2)\ne 0$. In this case, we see that $\Delta
M_Z^{(1,2)}=-1\ne \Delta M_X^{(0,2)}-1=0$.
\item Let us consider a scheme $X$ of degree $4$ with $2$ points on a
  $(1,0)$-line $R_0$ and other $2$ points on a $(0,1)$-line $C_0$ and its first
  difference $\Delta M_X$.
  \begin{figure}[H]
\begin{preview}
\begin{center}
\subfloat[$X$]{
\begin{tikzpicture}[line cap=round,line join=round,>=triangle 45,x=0.5cm,y=0.5cm]
\clip(-1,-1) rectangle (4,5);
\fill [color=black] (1,1) circle (2pt);
\fill [color=black] (1,2) circle (2pt);
\fill [color=black] (2,3) circle (2pt);
\fill [color=black] (3,3) circle (2pt);
\draw (1,0.5) -- (1,3.5);
\draw (2,2.5) -- (2,3.5);
\draw (3,2.5) -- (3,3.5);
\draw (0.5,1) -- (1.5,1);
\draw (0.5,2) -- (1.5,2);
\draw (0.5,3) -- (3.5,3);
\draw[color=black] (1,4.5) node {\footnotesize $C_0$};
\draw[color=black] (2,4.5) node {\footnotesize $C_1$};
\draw[color=black] (3,4.5) node {\footnotesize $C_2$};
\draw[color=black] (-0.5,1) node {\footnotesize $R_2$};
\draw[color=black] (-0.5,2) node {\footnotesize $R_1$};
\draw[color=black] (-0.5,3) node {\footnotesize $R_0$};
\end{tikzpicture}
}
\hspace{1cm}
\subfloat[$\Delta M_X$]{
\begin{tikzpicture}[x=0.7cm,y=0.7cm]
\clip(0,0) rectangle (7,7);
  \draw[style=help lines,xstep=1,ystep=1] (1,1) grid (6,6);
\foreach \x in {1,...,5} \draw (\x,1) +(.5,.5)  node {\dots};
\foreach \y in {2,...,5} \draw (5,\y) +(.5,.5) node {\dots};
\foreach \x in {1,...,4} \draw (\x,2.5) +(.5,0) node {$0$};
\foreach \y in {3,...,5} \draw (4.5,\y) +(0,.5) node {$0$};
\draw (1.5,3.5) node {$1$};
\draw (2.5,3.5) node {$0$};
\draw (3.5,3.5) node {$-1$};
\draw (1.5,4.5) node {$1$};
\draw (2.5,4.5) node {$0$};
\draw (3.5,4.5) node {$0$};
\foreach \x in {1,...,3} \draw (\x,5.5) +(.5,0) node {$1$};
\foreach \x in {0,...,4} \draw (\x,6.5) +(1.5,0) node {$\x$};
\foreach \y in {0,...,4} \draw (0.5,5.5-\y) node {$\y$};
\end{tikzpicture}
}
\end{center}
\end{preview}
\end{figure}
If $R$ is a $(1,0)$-line disjoint from $X$ and $Z=X\cup (R\cap C_2)$, then, under the notation of Theorem
\ref{T:1}, we have $r=2$, $n=2$, $q_1=1$, $q_2=2$ and $\Delta
M_X^{(2,2)}\ne 0$. However, $\Delta M_Z^{(2,1)}=0\ne \Delta M_X^{(1,1)}-1=-1$.
\begin{figure}[H]
\begin{preview}
\begin{center}
\subfloat[$Z$]{
\begin{tikzpicture}[line cap=round,line join=round,>=triangle 45,x=0.5cm,y=0.5cm]
\clip(-1,-1) rectangle (4,6);
\fill [color=black] (1,1) circle (2pt);
\fill [color=black] (1,2) circle (2pt);
\fill [color=black] (2,3) circle (2pt);
\fill [color=black] (3,3) circle (2pt);
\fill [color=black] (3,4) circle (2pt);
\draw (1,0.5) -- (1,4.5);
\draw (2,2.5) -- (2,4.5);
\draw (3,2.5) -- (3,4.5);
\draw (0.5,1) -- (1.5,1);
\draw (0.5,2) -- (1.5,2);
\draw (0.5,3) -- (3.5,3);
\draw (0.5,4) -- (3.5,4);
\draw[color=black] (1,5.5) node {\footnotesize $C_0$};
\draw[color=black] (2,5.5) node {\footnotesize $C_1$};
\draw[color=black] (3,5.5) node {\footnotesize $C_2$};
\draw[color=black] (-0.5,1) node {\footnotesize $R_2$};
\draw[color=black] (-0.5,2) node {\footnotesize $R_1$};
\draw[color=black] (-0.5,3) node {\footnotesize $R_0$};
\draw[color=black] (-0.5,4) node {\footnotesize $R$};
\end{tikzpicture}
}
\hspace{1cm}
\subfloat[$\Delta M_Z$]{
\begin{tikzpicture}[x=0.7cm,y=0.7cm]
\clip(0,0) rectangle (7,8);
  \draw[style=help lines,xstep=1,ystep=1] (1,1) grid (6,7);
\foreach \x in {1,...,5} \draw (\x,1) +(.5,.5)  node {\dots};
\foreach \y in {2,...,6} \draw (5,\y) +(.5,.5) node {\dots};
\foreach \x in {1,...,4} \draw (\x,2.5) +(.5,0) node {$0$};
\foreach \y in {3,...,6} \draw (4.5,\y) +(0,.5) node {$0$};
\draw (1.5,3.5) node {$1$};
\draw (2.5,3.5) node {$-1$};
\draw (3.5,3.5) node {$0$};
\draw (1.5,4.5) node {$1$};
\draw (2.5,4.5) node {$0$};
\draw (3.5,4.5) node {$-1$};
\foreach \x in {1,...,2} \draw (\x,5.5) +(.5,0) node {$1$};
\draw (3.5,5.5) node {$0$};
\foreach \x in {1,...,3} \draw (\x,6.5) +(.5,0) node {$1$};
\foreach \x in {0,...,4} \draw (\x,7.5) +(1.5,0) node {$\x$};
\foreach \y in {0,...,5} \draw (0.5,6.5-\y) node {$\y$};
\end{tikzpicture}
}
\end{center}
\end{preview}
\end{figure}  
\end{enumerate}
\end{ex}

\section{ACM case}

In this section we show that, if $X$ is an ACM scheme, then Theorem \ref{T:1} and Theorem \ref{T:2}
hold without any further assumption on $X$. The following result is
well known, but it is difficult to find a good reference and so we
give a short proof here.

\begin{prop}  \label{P}
  Let $X$ be an ACM $0$-dimensional scheme. Let $p_i=\#(X\cap R_i)$, for
  $i=0$,\dots,$a$  and let $q_j=\#(X\cap C_j)$, for
  $j=0$,\dots,$b$. Then:
\[
\Delta M_X^{(i,j)}=
\begin{cases}
  1 & \text{if } i\le q_j-1\text{ and }0\le j\le b\\
 0 & \text{otherwise}
\end{cases}
\]
or equivalently:
\[
\Delta M_X^{(i,j)}=
\begin{cases}
  1 & \text{if } j\le p_i-1\text{ and }0\le i\le a\\
 0 & \text{otherwise.}
\end{cases}
\]
\end{prop}
\begin{proof}
We show that:
\[
\Delta M_X^{(i,j)}=
\begin{cases}
  1 & \text{if } i\le q_j-1\text{ and }0\le j\le b\\
 0 & \text{otherwise}
\end{cases}
\]
The proof that also:
\[
\Delta M_X^{(i,j)}=
\begin{cases}
  1 & \text{if } j\le p_i-1\text{ and }0\le i\le a\\
 0 & \text{otherwise}
\end{cases}
\]
is similar.

  It is well known (see, for example, \cite{GM}, \cite{Gu} and
  \cite{VT}) that $X$ can be described after a suitable permutation of
  lines in such a way that the following conditions holds:
  \begin{enumerate}
  \item for every $i\in \{0,\dots,a\}$ there exists $j(i)\in
    \{0,\dots,b\}$ such that $R_i\cap C_j\in X$ for $j\in
    \{0,\dots,j(i)\}$ and $R_i\cap C_j\notin X$ for $j\in
    \{j(i)+1,\dots,b\}$;
\item $j(0)\ge j(1)\ge \dots \ge j(a)$.
  \end{enumerate}
Moreover, if any scheme $X$ satisfies these conditions, then $X$ is an
ACM scheme. Using this fact we can easily compute $\Delta M_X$ by induction on $a$. If
$a=0$, then the equality follows by the fact that $h^0(\mathscr
I_X(0,b))=0$ and $h^0(\mathscr I_X(0,b+1))=1$. 

If the equality holds for $a-1$, then we apply Theorem \ref{T2} and we
get the equality. 
\end{proof}

\begin{thm}  \label{T:5}
Let $X$ be an ACM scheme and let $R$ be a $(1,0)$-line disjoint from $X$. Let $C_{b+1}$,\dots,$C_n$,
$n\ge b$, be arbitrary $(0,1)$-lines and $i_1$,\dots,$i_r\in
\{0,\dots,b\}$. Let $\mathcal P=\{R\cap C_i\mid i\in\{0,\dots,n\},\,
i\ne i_1,\dots,i_r\}$ and let $Z=X\cup \mathcal P$. Suppose also that on the $(0,1)$-line $C_{i_k}$
there are $q_k$ points of $X$ for $k=1,\dots,r$ and that $q_1\le
q_2\le \dots \le q_r$. Then, given
$T=\{(q_1,n),(q_2,n-1),\dots,(q_r,n-r+1)\}$, we have:
\[
\Delta M_Z^{(i,j)}=
\begin{cases}
1 & \text{if } i=0,\, j\le n\\
0 & \text{if }i=0,\, j\ge n+1\\
 \Delta M_X^{(i-1,j)} & \text{if } i\ge 1\text{ and }(i,j)\notin T\\
 \Delta M_X^{(i-1,j)}-1 & \text{if } i\ge 1\text{ and }(i,j)\in T.
\end{cases}
\]
\end{thm}
\begin{proof}
The conclusion follows by Theorem \ref{T:1}, by Proposition \ref{P} and by
the fact that $\#(C_{b-k+1}\cap X)\le q_k$.
\end{proof}

In the same way:
\begin{thm}  \label{T:6}
Let $X$ be an ACM scheme and let $C$ be a $(0,1)$-line disjoint from $X$. Let $R_{a+1}$,\dots,$R_m$,
$m\ge a$, be arbitrary $(1,0)$-lines and $j_1$,\dots,$j_r\in
\{0,\dots,a\}$. Let $\mathcal P=\{C\cap R_j\mid j\in\{0,\dots,m\},\,
j\ne j_1,\dots,j_r\}$ and let $Z=X\cup \mathcal P$. Suppose also that on the $(1,0)$-line $R_{j_k}$
there are $p_k$ points of $X$ for $k=1,\dots,r$ and that $p_1\le
p_2\le \dots \le p_r$. Then, given
$T=\{(m,p_1),(m-1,p_2),\dots,(m-r+1,p_r)\}$, we have:
\[
\Delta M_Z^{(i,j)}=
\begin{cases}
1 & \text{if } i\le m,\, j=0\\
0 & \text{if }i\ge m+1,\, j=0\\
 \Delta M_X^{(i,j-1)} & \text{if } j\ge 1\text{ and }(i,j)\notin T\\
 \Delta M_X^{(i,j-1)}-1 & \text{if } j\ge 1\text{ and }(i,j)\in T.
\end{cases}
\]
\end{thm}
\begin{proof}
The works as in Theorem \ref{T:3}.  
\end{proof}

\section{Example}  \label{Esempio}

Now we show how it is possible to apply Theorem \ref{T:1}, Theorem
\ref{T:3} and Theorem \ref{T:5} to compute the first difference of the
Hilbert matrix of a scheme $X$ whose points can be distributed on a
grid of $(1,0)$ and $(0,1)$-lines in the following way:
\begin{figure}[H]
\begin{preview}
\begin{tikzpicture}[line cap=round,line join=round,>=triangle 45,x=0.5cm,y=0.5cm]
\foreach \x in {-1,1,2,3,4,5,6,7,8,9}
\foreach \y in {-1,1,2,3,4,5,6,7,8}
\clip(-1,0) rectangle (10,10);
\fill [color=black] (1,1) circle (2pt);
\fill [color=black] (2,1) circle (2pt);
\fill [color=black] (3,1) circle (2pt);
\fill [color=black] (1,2) circle (2pt);
\fill [color=black] (2,2) circle (2pt);
\fill [color=black] (3,2) circle (2pt);
\fill [color=black] (4,2) circle (2pt);
\fill [color=black] (1,3) circle (2pt);
\fill [color=black] (2,3) circle (2pt);
\fill [color=black] (3,3) circle (2pt);
\fill [color=black] (4,3) circle (2pt);
\fill [color=black] (5,3) circle (2pt);
\fill [color=black] (1,4) circle (2pt);
\fill [color=black] (2,4) circle (2pt);
\fill [color=black] (3,4) circle (2pt);
\fill [color=black] (4,4) circle (2pt);
\fill [color=black] (5,4) circle (2pt);
\fill [color=black] (6,5) circle (2pt);
\fill [color=black] (1,6) circle (2pt);
\fill [color=black] (2,6) circle (2pt);
\fill [color=black] (5,6) circle (2pt);
\fill [color=black] (6,6) circle (2pt);
\fill [color=black] (1,7) circle (2pt);
\fill [color=black] (2,7) circle (2pt);
\fill [color=black] (4,7) circle (2pt);
\fill [color=black] (5,7) circle (2pt);
\fill [color=black] (6,7) circle (2pt);
\fill [color=black] (3,8) circle (2pt);
\fill [color=black] (7,8) circle (2pt);
\fill [color=black] (8,8) circle (2pt);
\fill [color=black] (9,8) circle (2pt);
\draw (1,0.5) -- (1,8.5);
\draw[color=black] (1,9.5) node {\footnotesize $C_0$};
\draw (2,0.5) -- (2,8.5);
\draw[color=black] (2,9.5) node {\footnotesize $C_1$};
\draw (3,0.5) -- (3,8.5);
\draw[color=black] (3,9.5) node {\footnotesize $C_2$};
\draw (4,1.5) -- (4,8.5);
\draw[color=black] (4,9.5) node {\footnotesize $C_3$};
\draw (5,2.5) -- (5,8.5);
\draw[color=black] (5,9.5) node {\footnotesize $C_4$};
\draw (6,4.5) -- (6,8.5);
\draw[color=black] (6,9.5) node {\footnotesize $C_5$};
\draw (7,7.5) -- (7,8.5);
\draw[color=black] (7,9.5) node {\footnotesize $C_6$};
\draw (8,7.5) -- (8,8.5);
\draw[color=black] (8,9.5) node {\footnotesize $C_7$};
\draw (9,7.5) -- (9,8.5);
\draw[color=black] (9,9.5) node {\footnotesize $C_8$};
\draw (0.5,1) -- (3.5,1);
\draw[color=black] (-0.5,1) node {\footnotesize $R_7$};
\draw (0.5,2) -- (4.5,2);
\draw[color=black] (-0.5,2) node {\footnotesize $R_6$};
\draw (0.5,3) -- (5.5,3);
\draw[color=black] (-0.5,3) node {\footnotesize $R_5$};
\draw (0.5,4) -- (5.5,4);
\draw[color=black] (-0.5,4) node {\footnotesize $R_4$};
\draw (0.5,5) -- (6.5,5);
\draw[color=black] (-0.5,5) node {\footnotesize $R_3$};
\draw (0.5,6) -- (6.5,6);
\draw[color=black] (-0.5,6) node {\footnotesize $R_2$};
\draw (0.5,7) -- (6.5,7);
\draw[color=black] (-0.5,7) node {\footnotesize $R_1$};
\draw (0.5,8) -- (9.5,8);
\draw[color=black] (-0.5,8) node {\footnotesize $R_0$};
\end{tikzpicture}
\caption{The scheme $X$}
\end{preview}
\end{figure}
We compute $\Delta M_X$ by adding the points of the $(1,0)$-lines. The points on $R_4$,
$R_5$, $R_6$ and $R_7$ are an aCM scheme,
so that, by using Proposition \ref{P}, we get its first difference:
\begin{figure}[H]
\begin{preview}
\begin{center}
\subfloat
{
\begin{tikzpicture}[line cap=round,line join=round,>=triangle 45,x=0.5cm,y=0.5cm]
\clip(-1,-1) rectangle (6,5);
\foreach \x in {1,...,3} \fill [color=black] (\x,1) circle (2pt); 
\foreach \x in {1,...,4} \fill [color=black] (\x,2) circle (2pt); 
\foreach \x in {1,...,5} \fill [color=black] (\x,3) circle (2pt); 
\foreach \x in {1,...,5} \fill [color=black] (\x,4) circle (2pt); 
\draw (1,0.5) -- (1,4.5);
\draw (2,0.5) -- (2,4.5);
\draw (3,0.5) -- (3,4.5);
\draw (4,1.5) -- (4,4.5);
\draw (5,2.5) -- (5,4.5);
\draw (0.5,1) -- (3.5,1);
\draw (0.5,2) -- (4.5,2);
\draw (0.5,3) -- (5.5,3);
\draw (0.5,4) -- (5.5,4);
\end{tikzpicture}
}
\hspace{1cm}
\subfloat{
\begin{tikzpicture}[x=0.7cm,y=0.7cm]
\clip(0,0) rectangle (9,8);
  \draw[style=help lines,xstep=1,ystep=1] (1,1) grid (8,7);
\foreach \x in {1,...,7} \draw (\x,1) +(.5,.5)  node {\dots};
\foreach \y in {2,...,6} \draw (7,\y) +(.5,.5) node {\dots};
\foreach \x in {1,...,6} \draw (\x,2.5) +(.5,0) node {$0$};
\foreach \y in {3,...,6} \draw (6.5,\y) +(0,.5) node {$0$};
\foreach \x in {1,...,3} \draw (\x,3.5) +(.5,0) node {$1$};
\foreach \x in {1,...,4} \draw (\x,4.5) +(.5,0) node {$1$};
\foreach \x in {1,...,5} \draw (\x,5.5) +(.5,0) node {$1$};
\foreach \x in {1,...,5} \draw (\x,6.5) +(.5,0) node {$1$};
\draw (4.5,3.5) node {$0$};
\draw (5.5,3.5) node {$0$};
\draw (5.5,4.5) node {$0$};
\foreach \x in {0,...,6} \draw (\x,7.5) +(1.5,0) node {$\x$};
\foreach \y in {0,...,5} \draw (0.5,6.5-\y) node {$\y$};
\end{tikzpicture}
}
\end{center}
\end{preview}
\end{figure}
Now we add the point on the line $R_3$ and by Theorem \ref{T:5} we
compute its first difference:
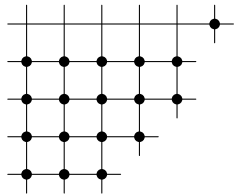
\begin{figure}[H]
\begin{preview}
\begin{center}
\subfloat
{
\begin{tikzpicture}[line cap=round,line join=round,>=triangle 45,x=0.5cm,y=0.5cm]
\clip(-1,-1) rectangle (7,6);
\foreach \x in {1,...,3} \fill [color=black] (\x,1) circle (2pt); 
\foreach \x in {1,...,4} \fill [color=black] (\x,2) circle (2pt); 
\foreach \x in {1,...,5} \fill [color=black] (\x,3) circle (2pt); 
\foreach \x in {1,...,5} \fill [color=black] (\x,4) circle (2pt); 
\fill [color=black] (6,5) circle (2pt); 
\draw (1,0.5) -- (1,5.5);
\draw (2,0.5) -- (2,5.5);
\draw (3,0.5) -- (3,5.5);
\draw (4,1.5) -- (4,5.5);
\draw (5,2.5) -- (5,5.5);
\draw (6,4.5) -- (6,5.5);
\draw (0.5,1) -- (3.5,1);
\draw (0.5,2) -- (4.5,2);
\draw (0.5,3) -- (5.5,3);
\draw (0.5,4) -- (5.5,4);
\draw (0.5,5) -- (6.5,5);
\end{tikzpicture}
}
\hspace{1cm}
\subfloat
{
\begin{tikzpicture}[x=0.7cm,y=0.7cm]
\clip(0,0) rectangle (10,9);
  \draw[style=help lines,xstep=1,ystep=1] (1,1) grid (9,8);
\foreach \x in {1,...,8} \draw (\x,1) +(.5,.5)  node {\dots};
\foreach \y in {2,...,7} \draw (8,\y) +(.5,.5) node {\dots};
\foreach \x in {1,...,7} \draw (\x,2.5) +(.5,0) node {$0$};
\foreach \y in {3,...,7} \draw (7.5,\y) +(0,.5) node {$0$};
\draw (1.5,3.5) node {$1$};
\draw (2.5,3.5) node {$0$};
\draw (3.5,3.5) node {$0$};
\draw (4.5,3.5) node {$-1$};
\draw (5.5,3.5) node {$0$};
\draw (6.5,3.5) node {$0$};
\foreach \x in {1,...,4} \draw (\x,4.5) +(.5,0) node {$1$};
\draw (5.5,4.5) node {$-1$};
\draw (6.5,4.5) node {$0$};
\foreach \x in {1,...,5} \draw (\x,5.5) +(.5,0) node {$1$};
\draw (6.5,5.5) node {$-1$};
\foreach \x in {1,...,5} \draw (\x,6.5) +(.5,0) node {$1$};
\draw (6.5,6.5) node {$0$};
\foreach \x in {1,...,6} \draw (\x,7.5) +(.5,0) node {$1$};
\foreach \x in {0,...,7} \draw (\x,8.5) +(1.5,0) node {$\x$};
\foreach \y in {0,...,6} \draw (0.5,7.5-\y) node {$\y$};
\end{tikzpicture}
}
\end{center}
\end{preview}
\end{figure}
In the same way we add the points on $R_2$ and by Theorem \ref{T:1} we
compute the first difference:
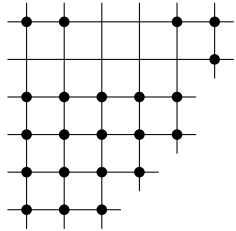
\begin{figure}[H]
\begin{preview}
\begin{center}
\subfloat
{
\begin{tikzpicture}[line cap=round,line join=round,>=triangle 45,x=0.5cm,y=0.5cm]
\clip(-1,-1) rectangle (7,7);
\foreach \x in {1,...,3} \fill [color=black] (\x,1) circle (2pt); 
\foreach \x in {1,...,4} \fill [color=black] (\x,2) circle (2pt); 
\foreach \x in {1,...,5} \fill [color=black] (\x,3) circle (2pt); 
\foreach \x in {1,...,5} \fill [color=black] (\x,4) circle (2pt); 
\fill [color=black] (6,5) circle (2pt); 
\foreach \x in {5,6} \fill [color=black] (\x,6) circle (2pt);
\foreach \x in {1,2} \fill [color=black] (\x,6) circle (2pt);
\draw (1,0.5) -- (1,6.5);
\draw (2,0.5) -- (2,6.5);
\draw (3,0.5) -- (3,6.5);
\draw (4,1.5) -- (4,6.5);
\draw (5,2.5) -- (5,6.5);
\draw (6,4.5) -- (6,6.5);
\draw (0.5,1) -- (3.5,1);
\draw (0.5,2) -- (4.5,2);
\draw (0.5,3) -- (5.5,3);
\draw (0.5,4) -- (5.5,4);
\draw (0.5,5) -- (6.5,5);
\draw (0.5,6) -- (6.5,6);
\end{tikzpicture}
}
\hspace{1cm}
\subfloat
{
\begin{tikzpicture}[x=0.7cm,y=0.7cm]
\clip(0,0) rectangle (10,10);
  \draw[style=help lines,xstep=1,ystep=1] (1,1) grid (9,9);
\foreach \x in {1,...,8} \draw (\x,1) +(.5,.5)  node {\dots};
\foreach \y in {2,...,8} \draw (8,\y) +(.5,.5) node {\dots};
\foreach \x in {1,...,7} \draw (\x,2.5) +(.5,0) node {$0$};
\foreach \y in {3,...,8} \draw (7.5,\y) +(0,.5) node {$0$};
\draw (1.5,3.5) node {$1$};
\draw (2.5,3.5) node {$0$};
\draw (3.5,3.5) node {$0$};
\draw (4.5,3.5) node {$-1$};
\draw (5.5,3.5) node {$0$};
\draw (6.5,3.5) node {$0$};
\foreach \x in {1,2} \draw (\x,4.5) +(.5,0) node {$1$};
\foreach \x in {3,4} \draw (\x,4.5) +(.5,0) node {$1$};
\draw (5.5,4.5) node {$-2$};
\draw (6.5,4.5) node {$0$};
\foreach \x in {1,...,5} \draw (\x,5.5) +(.5,0) node {$1$};
\draw (6.5,5.5) node {$-2$};
\foreach \x in {1,...,5} \draw (\x,6.5) +(.5,0) node {$1$};
\draw (6.5,6.5) node {$0$};
\foreach \x in {1,...,6} \draw (\x,7.5) +(.5,0) node {$1$};
\foreach \x in {1,...,6} \draw (\x,8.5) +(.5,0) node {$1$};
\foreach \x in {0,...,7} \draw (\x,9.5) +(1.5,0) node {$\x$};
\foreach \y in {0,...,7} \draw (0.5,8.5-\y) node {$\y$};
\end{tikzpicture}
}
\end{center}
\end{preview}
\end{figure}
Now we add the points on $R_1$ and again by Theorem \ref{T:1} we
compute its first difference:

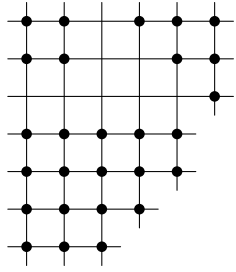
\begin{figure}[H]
\begin{preview}
\begin{center}
\subfloat
{
\begin{tikzpicture}[line cap=round,line join=round,>=triangle 45,x=0.5cm,y=0.5cm]
\clip(-1,-1) rectangle (7,8);
\foreach \x in {1,...,3} \fill [color=black] (\x,1) circle (2pt); 
\foreach \x in {1,...,4} \fill [color=black] (\x,2) circle (2pt); 
\foreach \x in {1,...,5} \fill [color=black] (\x,3) circle (2pt); 
\foreach \x in {1,...,5} \fill [color=black] (\x,4) circle (2pt); 
\fill [color=black] (6,5) circle (2pt); 
\foreach \x in {1,2} \fill [color=black] (\x,6) circle (2pt);
\foreach \x in {5,6} \fill [color=black] (\x,6) circle (2pt);
\fill [color=black] (1,7) circle (2pt);
\fill [color=black] (2,7) circle (2pt);
\foreach \x in {4,...,6} \fill [color=black] (\x,7) circle (2pt);
\draw (1,0.5) -- (1,7.5);
\draw (2,0.5) -- (2,7.5);
\draw (3,0.5) -- (3,7.5);
\draw (4,1.5) -- (4,7.5);
\draw (5,2.5) -- (5,7.5);
\draw (6,4.5) -- (6,7.5);
\draw (0.5,1) -- (3.5,1);
\draw (0.5,2) -- (4.5,2);
\draw (0.5,3) -- (5.5,3);
\draw (0.5,4) -- (5.5,4);
\draw (0.5,5) -- (6.5,5);
\draw (0.5,6) -- (6.5,6);
\draw (0.5,7) -- (6.5,7);
\end{tikzpicture}
}
\hspace{1cm}
\subfloat
{
\begin{tikzpicture}[x=0.7cm,y=0.7cm]
\clip(0,0) rectangle (10,11);
  \draw[style=help lines,xstep=1,ystep=1] (1,1) grid (9,10);
\foreach \x in {1,...,8} \draw (\x,1) +(.5,.5)  node {\dots};
\foreach \y in {2,...,9} \draw (8,\y) +(.5,.5) node {\dots};
\foreach \x in {1,...,7} \draw (\x,2.5) +(.5,0) node {$0$};
\foreach \y in {3,...,9} \draw (7.5,\y) +(0,.5) node {$0$};
\draw (1.5,3.5) node {$1$};
\draw (2.5,3.5) node {$0$};
\draw (3.5,3.5) node {$0$};
\draw (4.5,3.5) node {$-1$};
\draw (5.5,3.5) node {$0$};
\draw (6.5,3.5) node {$0$};
\foreach \x in {1,2} \draw (\x,4.5) +(.5,0) node {$1$};
\foreach \x in {3,4} \draw (\x,4.5) +(.5,0) node {$1$};
\draw (5.5,4.5) node {$-2$};
\draw (6.5,4.5) node {$0$};
\foreach \x in {1,...,5} \draw (\x,5.5) +(.5,0) node {$1$};
\draw (6.5,5.5) node {$-3$};
\foreach \x in {1,...,5} \draw (\x,6.5) +(.5,0) node {$1$};
\draw (6.5,6.5) node {$0$};
\foreach \x in {1,...,6} \draw (\x,7.5) +(.5,0) node {$1$};
\foreach \x in {1,...,6} \draw (\x,8.5) +(.5,0) node {$1$};
\foreach \x in {1,...,6} \draw (\x,9.5) +(.5,0) node {$1$};
\foreach \x in {0,...,7} \draw (\x,10.5) +(1.5,0) node {$\x$};
\foreach \y in {0,...,8} \draw (0.5,9.5-\y) node {$\y$};
\end{tikzpicture}
}
\end{center}
\end{preview}
\end{figure}

Finally, by applying again Theorem \ref{T:1} we get the first
difference $\Delta M_X$ of $X$.
\begin{figure}[H]
\begin{preview}
\begin{center}
\subfloat{
\begin{tikzpicture}[line cap=round,line join=round,>=triangle 45,x=0.5cm,y=0.5cm]
\foreach \x in {-1,1,2,3,4,5,6,7,8,9}
\foreach \y in {-1,1,2,3,4,5,6,7,8}
\clip(-1,0) rectangle (10,10);
\fill [color=black] (1,1) circle (2pt);
\fill [color=black] (2,1) circle (2pt);
\fill [color=black] (3,1) circle (2pt);
\fill [color=black] (1,2) circle (2pt);
\fill [color=black] (2,2) circle (2pt);
\fill [color=black] (3,2) circle (2pt);
\fill [color=black] (4,2) circle (2pt);
\fill [color=black] (1,3) circle (2pt);
\fill [color=black] (2,3) circle (2pt);
\fill [color=black] (3,3) circle (2pt);
\fill [color=black] (4,3) circle (2pt);
\fill [color=black] (5,3) circle (2pt);
\fill [color=black] (1,4) circle (2pt);
\fill [color=black] (2,4) circle (2pt);
\fill [color=black] (3,4) circle (2pt);
\fill [color=black] (4,4) circle (2pt);
\fill [color=black] (5,4) circle (2pt);
\fill [color=black] (6,5) circle (2pt);
\fill [color=black] (1,6) circle (2pt);
\fill [color=black] (2,6) circle (2pt);
\fill [color=black] (5,6) circle (2pt);
\fill [color=black] (6,6) circle (2pt);
\fill [color=black] (1,7) circle (2pt);
\fill [color=black] (2,7) circle (2pt);
\fill [color=black] (4,7) circle (2pt);
\fill [color=black] (5,7) circle (2pt);
\fill [color=black] (6,7) circle (2pt);
\fill [color=black] (3,8) circle (2pt);
\fill [color=black] (7,8) circle (2pt);
\fill [color=black] (8,8) circle (2pt);
\fill [color=black] (9,8) circle (2pt);
\draw (1,0.5) -- (1,8.5);
\draw (2,0.5) -- (2,8.5);
\draw (3,0.5) -- (3,8.5);
\draw (4,1.5) -- (4,8.5);
\draw (5,2.5) -- (5,8.5);
\draw (6,4.5) -- (6,8.5);
\draw (7,7.5) -- (7,8.5);
\draw (8,7.5) -- (8,8.5);
\draw (9,7.5) -- (9,8.5);
\draw (0.5,1) -- (3.5,1);
\draw (0.5,2) -- (4.5,2);
\draw (0.5,3) -- (5.5,3);
\draw (0.5,4) -- (5.5,4);
\draw (0.5,5) -- (6.5,5);
\draw (0.5,6) -- (6.5,6);
\draw (0.5,7) -- (6.5,7);
\draw (0.5,8) -- (9.5,8);
\end{tikzpicture}
}
\hspace{1cm}
\subfloat
{
\begin{tikzpicture}[x=0.7cm,y=0.7cm]
\clip(0,0) rectangle (13,12);
  \draw[style=help lines,xstep=1,ystep=1] (1,1) grid (12,11);
\foreach \x in {1,...,11} \draw (\x,1) +(.5,.5)  node {\dots};
\foreach \y in {2,...,10} \draw (11,\y) +(.5,.5) node {\dots};
\foreach \x in {1,...,10} \draw (\x,2.5) +(.5,0) node {$0$};
\foreach \y in {3,...,10} \draw (10.5,\y) +(0,.5) node {$0$};
\draw (1.5,3.5) node {$1$};
\draw (2.5,3.5) node {$0$};
\draw (3.5,3.5) node {$0$};
\draw (4.5,3.5) node {$-1$};
\foreach \x in {5,...,9} \draw (\x,3.5) +(.5,0) node {$0$};
\foreach \x in {1,2} \draw (\x,4.5) +(.5,0) node {$1$};
\foreach \x in {3,4} \draw (\x,4.5) +(.5,0) node {$1$};
\draw (5.5,4.5) node {$-3$};
\draw (6.5,4.5) node {$-1$};
\foreach \x in {7,...,9} \draw (\x,4.5) +(.5,0) node {$0$};
\foreach \x in {1,...,5} \draw (\x,5.5) +(.5,0) node {$1$};
\draw (6.5,5.5) node {$-3$};
\draw (7.5,5.5) node {$0$};
\draw (8.5,5.5) node {$0$};
\draw (9.5,5.5) node {$0$};
\foreach \x in {1,...,5} \draw (\x,6.5) +(.5,0) node {$1$};
\draw (6.5,6.5) node {$0$};
\draw (7.5,6.5) node {$-1$};
\draw (8.5,6.5) node {$-1$};
\draw (9.5,6.5) node {$0$};
\foreach \x in {1,...,6} \draw (\x,7.5) +(.5,0) node {$1$};
\foreach \x in {7,8} \draw (\x,7.5) +(.5,0) node {$0$};
\draw (9.5,7.5) node {$-1$};
\foreach \x in {1,...,6} \draw (\x,8.5) +(.5,0) node {$1$};
\foreach \x in {7,...,9} \draw (\x,8.5) +(.5,0) node {$0$};
\foreach \x in {1,...,6} \draw (\x,9.5) +(.5,0) node {$1$};
\foreach \x in {7,...,9} \draw (\x,9.5) +(.5,0) node {$0$};
\foreach \x in {1,...,9} \draw (\x,10.5) +(.5,0) node {$1$};
\foreach \x in {0,...,10} \draw (\x,11.5) +(1.5,0) node {$\x$};
\foreach \y in {0,...,9} \draw (0.5,10.5-\y) node {$\y$};
\end{tikzpicture}
}
\end{center}
\end{preview}
\end{figure}

\end{document}